\documentclass[11pt,leqno]{article}

\usepackage{latexsym}
\usepackage{amssymb,amsbsy,amsmath,amsfonts,amscd}
\usepackage{vmargin}
\usepackage{epsfig}
\usepackage{amsmath}
\usepackage{amsfonts}
\usepackage{amssymb}
\usepackage{graphicx}%
\usepackage[dvips]{color}

\newtheorem{lemme}{Lemma}[section]

\newtheorem{Theorem}{Theorem}[section]

\newtheorem{Proposition}{Proposition}[section]
\newtheorem{Remark}{Remark}[section]

\newtheorem{Corollary}{Corollary}[section]

\newcommand{ \vit}{{\bu}}

\newcommand{\R}{\ifmmode{{\rm I} \hskip -2pt {\rm R}}
    \else{\hbox{$I\hskip -2pt R$}}\fi}
\newcommand{\Z}{\mathbb{Z}}

\def\vec#1{\boldsymbol{#1}}
\def\bu{\vec{u}}
\def\Lim{\mathop{\mathrm{Lim}}\nolimits}
\begin{document}
\begin{titlepage}
\title{\vskip1cm\textbf{Ladder Theorem and length-scale estimates for a Leray alpha model of turbulence}\vskip0.1cm}
\author{\begin{Large}$\hbox{Hani Ali}\thanks{IRMAR, UMR 6625,
Universit\'e Rennes 1,
Campus Beaulieu,
35042 Rennes cedex
FRANCE;
hani.ali@univ-rennes1.fr}$ \end{Large}}
\end{titlepage}
\date{}
\maketitle
\begin{center}
 \vskip-0.1cm\textbf{Abstract}\vskip0.1cm
\end{center}
\  \ In this paper, we study the Modified Leray alpha model with  periodic boundary conditions. 
We show that the regular solution verifies a sequence of energy inequalities that is called "ladder inequalities".
 Furthermore, we estimate some quantities of physical relevance in terms of the Reynolds number.\\

\vskip0.05cm
\section{Introduction}
 We consider  in this paper the ML-$\alpha$ model of turbulence 
\begin{equation}
\label{Malpha ns}
 \left\{
\begin{array} {llll} \displaystyle
 \frac{\partial \vec{u}^{}}{\partial t}+({\vec{u}^{}} \cdot \nabla)
\overline{\vec{u}^{}} - \nu \Delta \vec{u}^{} + \nabla p^{} =
{\vec{f}^{}} \ \ \ \ \ \hbox{in}\ \R^{+}\times\mathbb{T}_3,\\
- \alpha^{2}\Delta\overline{\vit}+\overline{\vit}+\nabla \pi  = \vec{u} \ \ \ \ \ \hbox{in}\ \mathbb{T}_3, \\
 \nabla \cdot \vec{u}^{}= \nabla \cdot \overline{\vit}=0,\\
 \displaystyle \oint_{\mathbb{T}_3} \vec{u}^{}=\displaystyle \oint_{\mathbb{T}_3}\overline{\vit}=0,\\
\vec{u}^{}_{t=0}={\vec{u}^{in}}.
\end{array}\right.
\end{equation}
The boundary conditions are periodic boundary conditions. Therefore we consider these equations on the  three dimensional torus $\mathbb{T}_3=\left ( \R^3 / {\cal T}_3 \right) $ where ${\cal T}_{3} = 2 \pi \Z^{3} /L$ , $\textbf{x} \in \mathbb{T}_3, \ \hbox{and} \ t \in ]0,+\infty[.$
The unkowns are  the  velocity vector field $ \vec{u}$
and the scalar
pressure  p. The viscosity $\nu,$   the initial velocity vector field $\vec{u}^{in},$ and   the external
force $\vec{f}$ with $\nabla \cdot \vec{f} = 0$ are given. In this paper the force $\vec{f}$ does not depend on time.\\
This model has been first studied in \cite{ILT05}, where the authors prove the global existence and uniqueness of the solution.
They also prove  the existence of a global attractor $ \mathcal{A}$  to this model and they made estimates of the fractal dimension of this attractor in terms of Grashof number $Gr$.\\
The dimension of the attractor gives us some idea of the level of the complexity of the flow. 
 The relation between the number of determining modes, determining nodes and the evolution of volume elements of the attractors are discussed by Jones and Titi in \cite{JT93}.  
Temam also  interprets in his book \cite{RT88} the dimension of the attractors as the number of degrees of freedom of the flow. 

It is easily seen that when $\alpha = 0$, eqs. (\ref{Malpha ns}) reduce to the usual Navier Stokes equations for incompressible fluids.\\
Assuming that $ \vec{f} \in C^{\infty},$   any $C^{\infty}$ solution  to the Navier Stokes equations verifies formally what is called the ladder inequality \cite{DG95}. That means,
for any  $C^{\infty}$ solution $(\vec{u},p)$ to the (NSE),  the velocity part $\vec{u}$ satisfies the following relation between its higher derivatives,
\begin{equation} \label{2.27}
\begin{array}{cc}
 \displaystyle \frac{1}{2}\frac{d}{dt}{H_{N}} \le -\nu H_{N+1} + C_{N}
H_{N} \left\|\nabla\vec{u}\right\|_{\infty} + H_{N}^{1/2}
\textbf{$\Phi_{N}$}^{1/2},\\
\hbox{where } \displaystyle H_{N}=\int_{\mathbb{T}_3}\left|\nabla^{N} {\vec{u}}\right|^{2} d\textbf{x} \hbox{ and } \displaystyle \textbf{$\Phi_{N}$} =\int_{\mathbb{T}_3}|
\nabla^{N}\vec{f}|^{2}d\textbf{x}. 
\end{array}
\end{equation}
This differential inequalities are used first in \cite{DG95} to show the existence of a lower bound on the smallest scale in the flow. The same result is obtained  in \cite{DT95}   by a Gevrey Class estimates. 
Recently, the  ladder inequalities are used  to study the intermittency of solutions to the    Navier Stokes equations   see \cite{GD05}.
While the ladder inequalities to the Navier Stokes equations are based on the assumption that a solution exists, so that the higher order norms are finite, no such assumption is necessary in the case of alpha regularistaion where existence and uniqueness of a $C^{\infty}$ solution  are guaranteed.     
The ladder inequalties are generalized in \cite{GH06,GH08} to other equations based on the Navier stokes equations such as Navier Stokes-alpha model \cite{FDT02} and Leray alpha model \cite{CHOT05}.\\
 We aim to study in this paper ladder inequalities for model (\ref{Malpha ns}).
 In the whole paper, $\alpha >0$ is given and we assume that  the initial data is $C^{\infty}$. One of the  main results of this paper is: 
 \begin{Theorem}
\label{thhhh}  Assume  $\vec{f} \in C^{\infty}(\mathbb{T}_3)^{3} $ and ${\vec{u}}^{in} \in C^{\infty}(\mathbb{T}_3)^{3}$. Let $(\vec{u},p):=( \vec{u}^{\alpha},p^{\alpha})$ be the unique solution to problem (\ref{Malpha ns}).Then the velocity part $\vec{u}$ satisfies
the ladder inequalities,
\begin{equation}
\label{inegalite d'echelle}
\begin{array}{ll}
 \displaystyle \frac{1}{2}(\frac{d}{dt}{\overline{H_{N}}}+\alpha^{2}\frac{d}{dt}{\overline{H_{N+1}}}) \le \displaystyle  -\nu (\overline{H_{N+1}}+\alpha^{2}\overline{H_{N+2}})\\
  \hskip 3 cm \displaystyle + C_{N}   \left\|\nabla\overline{\vec{u}}\right\|_{\infty} (\overline{H_{N}} + \displaystyle \alpha^{2} \overline{H_{N+1}}) + \overline{{H_{N}}}^{1/2} \textbf{$\Phi_{N}$}^{1/2},
 \end{array}
\end{equation}
where
 \begin{equation}
\label{propriete10}
 \overline{H_{N}}=\int_{\mathbb{T}_3}\left|\nabla^{N} \overline{\vec{u}}\right|^{2} d\textbf{x}, C_0= 0 \hbox{ and } C_N \approx 2^N \hbox{ for all } N \ge 1.   
\end{equation}
\end{Theorem}
The gradient symbol $ \nabla^N$ here refers to all derivatives of evrey component of $\bu$ of order $N$ in $L^2(\mathbb{T}_3)$.
\begin{Remark}
We note that, as $\alpha  \rightarrow 0$,  $\overline{H_{N}} \rightarrow {H_{N}}$. Thus  we find the inequality (\ref{2.27}).
\end{Remark}
 
 Another Task of this paper is to estimate quantities of physical relevance in terms of the Reynolds number,  these result are summarized in the table \ref{Table1} whose proof is given in section 5.
 For simplicity the eqs. (\ref{Malpha ns}) will be considered with forcing $ \vec{f}(\vec{x})$ taken to be $L^2$ bounded of narrow band type with a single lenght scale $\ell$ (see \cite{GD05, GH06}) such that 
\begin{equation}
\|\nabla^n  \vec{f}  \|_{L^2} \approx \ell^{-n} \|  \vec{f}  \|_{L^2}.
\end{equation}
In order to estimate small length sacles associated with higher order moments, we  combine in section 5 the force with the higher derivative of the velocity such that 
 \begin{equation}
J_N = \overline{F_N} + 2 \alpha^2 \overline{F_{N+1}},
\end{equation}
where 
 \begin{equation} 
 \overline{F_N} = \overline{H_N} + \tau {\Phi_{N}},
  \end{equation}
and  the quantity $ \tau$ is defined by
 \begin{equation}
\tau = \ell^2 \nu^{-1} (Gr \ln{Gr} )^{-1/2}.
\end{equation}
The ${J_N} $  is used to define a set of time-depend inverse length scales  
  \begin{equation}
  \kappa_{N,r}=\left(  \frac{J_N}{J_r} \right)^{\frac{1}{2(N-r)}}. 
  \end{equation}

The second main result of the paper is the  following Theorem. 

\begin{Theorem}
\label{lengththoerem}
 Let $\vec{f} \in C^{\infty}(\mathbb{T}_3)^{3} $ with narrow-band type and ${\vec{u}}^{in} \in C^{\infty}(\mathbb{T}_3)^{3}$. Let $ \vec{u}:=\vec{u}^{\alpha}$ be the velocity part of the  solution to problem (\ref{Malpha ns}). Then  estimates in term of Reynolds number $Re$ for the length sacles associated with higher order moments solution $ \kappa_{N,0}$ $(N \ge  2 )$, the inverse Kolomogrov length $ \lambda_k$ and the attractor dimesion $ d_{F, ML-\alpha}(\mathcal{A})$  are given by 
    \begin{equation}
   \label{lengthimproved}
   \ell^{2} \left\langle   \kappa_{N,0}^2 \right\rangle 
\le   C(\alpha,\nu , \ell, L)^{(N-1)/N} Re^{5/2-3/2N} (\ln{Re})^{1/N}+ C Re \ln{Re}.
\end{equation}
 \begin{equation}
 \label{improvedbound kolomogrov}
 \ell \lambda_k^{-1} \le c Re^{5/8}.
 \end{equation}
 \begin{equation}
 \label{attracteurimproved}
d_{F, ML-\alpha}(\mathcal{A}) \le c \left( \frac{L^3 \ell^{-4}}{\alpha^2 \lambda_1^{3/2}}  \right)^{3/4} Re^{9/4}.
\end{equation}
Where $\left\langle \cdot \right\rangle $ is the long time average defined below (\ref{longtime})
\end{Theorem}

\ \ The paper is organized as follows: In section 2, we start by summarizing and discussing the results given above. In section 3 we recall some helpfuls results about  existence and uniqueness  for this
ML$-\alpha$ model,
and  we prove  a general regularity result. 
 In section 4, we prove Theorem \ref{thhhh}. We stress that for all $N \in \mathbb{N} $ fixed, inequality (\ref{inegalite d'echelle}) goes to inequality (\ref{2.27}) when  $\alpha \rightarrow 0,$ at least formally.
In section 5, we prove Theorem \ref{lengththoerem}.

%

\section{Summary and discussion of the results}
%
  
%

Generally the most important of the estimates in Navier stokes theory have been found in terms of the Grashof number $Gr$ defined below in terms of the forcing, but these are difficult to compare with the results of Kolomogrov scaling theories \cite{UF95} which are expressed in terms of Reynolds number $Re$ based on the Navier Stokes velocity $ \vec{u}$. A good definition of this is 
\begin{equation}
Re=\frac{U \ell}{\nu}, \ \ U^2 = L^{-3}\left\langle \|\vec{u} \|_{L^2}^2 \right\rangle,
\end{equation}
where $\left\langle \cdot \right\rangle $ is the long time average 
\begin{equation}
\label{longtime}
\left\langle g( \cdot) \right\rangle = {\Lim}_{t\rightarrow\infty}\frac{1}{t}\int_{0}^{t}g(s)ds.
\end{equation}
Where $\Lim$ indicates a generalized limit that extends the usual limits \cite{FMRT01}.\\
With $f_{rms}=L^{-3/2}\|\vec{f} \|_{L^2},$ the standard definition of the Grashof number in three dimensions is 
\begin{equation}
Gr= \frac{\ell^{3} f_{rms}}{\nu^2}.
\end{equation}
 Doering and Foias \cite{DF02} have addressed the problem of how to relate $ Re$ and $Gr$ and have shown that in the limit $Gr \rightarrow \infty$, solutions of the Navier Stokes equations must satisfy 
 \begin{equation}
 \label{gr and re}
 Gr \le c (Re^2 + Re).
 \end{equation}
 
 Using the above relation (\ref{gr and re}),  Doering and Gibbon \cite{GD05}  have re-expressed some Navier Stokes estimates in terms of $Re$. In particular they showed that  the energy dissipation rate  $ \epsilon=\nu \left\langle  \|\nabla \vec{u}\|_{L^2}^2\right\rangle L^{-3} $  is bounded above by 
  \begin{equation}
 \label{energy disss}
\epsilon  \le  c \nu^3 \ell^{-4} \left( Re^3 + Re \right),
 \end{equation}
 
  and the inverse kolomogrov length $ \lambda_k^{-1}= (\epsilon/\nu^3)^{1/4}$ is  bounded above by 
 \begin{equation}
 \label{bound kolomogrov}
 \ell \lambda_k^{-1} \le c Re^{3/4}.
 \end{equation}
 
 The relation (\ref{gr and re}) is essentially a Navier Stokes result.  In \cite{GH06}  it has been shown that this property holds for the Navier Stokes-alpha model \cite{FDT02}; the same methods can be used  to show this also holds for eqs. (\ref{Malpha ns}). In this paper, we will use  (\ref{gr and re})  to obtain estimates in terms of the Reynolds number $Re$.
 \begin{table}[h]
\begin{tabular}{|l|c|c|c|c|r|}
  \hline
   & \small{NS} & \small{NS-}$\alpha$/\small{Bardina}& \small{Leray}-$\alpha$  & \small{ML-}$\alpha$ & \small{Eq.} \\
  \hline
   \hline
$ \ell \lambda_k^{-1}$ & $Re^{3/4}$ &$ Re^{5/8}$ &$ Re^{7/12}$ &$ Re^{5/8}$ &(\ref{bound kolomogrov}) \\
 $\left\langle \overline{H_1} \right\rangle$ &  $Re^{3}$    & $Re^{5/2}$     & $Re^{7/3}$                      & $Re^{5/2}$ &  (\ref{H 1 bar})  \\
  $\left\langle \overline{H_2} \right\rangle$ &   -        & $Re^{3}$     & $Re^{8/3}$                & $Re^{3}$  &  (\ref{H 2 bar})   \\
   $\left\langle \overline{H_3} \right\rangle$ &     -            &    - / -       & $ Re^{3}$           & $ Re^{7}$&   (\ref{H 3 bar})   \\
    $d_F(\mathcal{A})$ &    -    & $Re^{9/4} \   / \  Re^{9/5}$           & $Re^{9/7}$           & $Re^{9/4}$ &  (\ref{attractor dimension}) \\ 
    $\ell^2\left\langle \kappa_{N,r}^2 \right\rangle $   & -      & $Re^{11/4} $        & $Re^{17/4} $                 & $Re^{5/2} $ &   (\ref{k n r}) \\
     $\ell^2\left\langle \kappa_{1,0}^2 \right\rangle $ &  $Re \ln{Re}$              &  $Re \ln{Re}$    &  $Re \ln{Re}$             &$Re \ln{Re}$ & (\ref{k 1 0})   \\
     $\left\langle \| \overline{\vec{u}}\|_{\infty}^2 \right\rangle $ & -    & $Re^{11/4} $     & $Re^{5/2} $            & $Re^{11/4} $&  (\ref{u infini})  \\
      $\left\langle \|\nabla \overline{\vec{u}}\|_{\infty} \right\rangle $ & - & $Re^{35/16}$   & $Re^{17/12}$  & $Re^{5/2}$&  (\ref{nabla u infini})\\
      $\ell^2 \left\langle \kappa_{N,0}^2 \right\rangle$&-&$Re^{\frac{11}{4}-\frac{7}{4N}} (\ln{Re})^{\frac{1}{N}}$&$Re^{\frac{17}{12}-\frac{5}{12N}}(\ln{Re})^{\frac{1}{N}}$             &$Re^{\frac{5}{2}-\frac{3}{2N}}(\ln{Re})^{\frac{1}{N}}$ &(\ref{N DEPEND})\\
     \hline
    \end{tabular}
    \caption{\small{Comparison of various upper bounds for the Navier Stokes, Navier Stokes-$\alpha$, Bardina, Leray-$\alpha$ and Modified Leray-$\alpha$ with constant omitted}}
    \label{Table1}
   \end{table}

  These estimates are listed in Table 1. The estimate for $d_{F, ML-\alpha}(\mathcal{A})$ are consistent with the long-standing belief that $ Re^{3/4}\times Re^{3/4}\times Re^{3/4}$ resolution grid points are needed to numerically resolve the flow.  The fact that this bound is not valid to the Navier Stokes equations is consistent with the fact that $d_{F, ML-\alpha}(\mathcal{A})$ blows up as $ \alpha$ tends to zero.  The improved estimate to the inverse kolomgrov $\lambda_k^{-1}$  coincide with the estimate to the Navier Stokes alpha given in \cite{GH06} and  blows up when $\alpha$ tends to zero.
 The estimate for $ \left\langle \kappa_{N,0}^2 \right\rangle $  comes out to be sharper  than those given for the Navier Stokes alpha because of the $\| \nabla \bu \|_{\infty}$ term in the ladder inequality as opposed to the $ \nu^{-1} \|\bu \|^{2}_{\infty}$ in   \cite{GH06}.
 This estimate gives us a length scale that is immensely small. Such scale is unreachable computationally and the regular solution on a neighbour of this scale is unresolvable.
 Thus the resolution issues in computations of the flow
are not only  associated with the problem of regularity but they also raise the
question of how resolution length scales can be defined and estimated.\\
   
 \ We finish this section by the following remark. The existence and the uniqueness of a $C^{\infty}$ solution for all time $T$ to the ML-$\alpha$  motivate the present study. Provided that regular solution exists for a  maximal inteval time $[0,T^{*}[$, we can show the ladder inequalities to the Navier Stokes equations in $[0,T^{*}[$.
We then naturally ask ourselves If can we use the convergence of  (\ref{inegalite d'echelle}) to  (\ref{2.27}) in   $[0,T^{*}[$    to deduce some informations about the regular solution beyond the time $ T^{*}$? This is  an crucial problem.


\section{ Existence, unicity  and Regularity results}
 We begin this section by recalling   the system (\ref{Malpha ns}) considered with periodic boundary conditions.
 \begin{equation}
\label{Malphaexpanded ns}
 \left\{
\begin{array} {llll} \displaystyle
 \frac{\partial \vec{u}^{}}{\partial t}+({\vec{u}^{}} \cdot \nabla)
\overline{\vec{u}^{}} - \nu \Delta \vec{u}^{} + \nabla p^{} =
{\vec{f}^{}} \ \ \ \ \ \hbox{in}\ \R^{+}\times\mathbb{T}_3,\\
- \alpha^{2}\Delta\overline{\vit}+\overline{\vit}+\nabla \pi  = \vec{u} \ \ \ \ \ \hbox{in}\ \mathbb{T}_3, \\
 \nabla \cdot \vec{u}^{}= \nabla \cdot \overline{\vit}=0,\\
 \displaystyle \oint_{\mathbb{T}_3} \vec{u}^{}=\displaystyle \oint_{\mathbb{T}_3}\overline{\vit}=0,\\
\vec{u}^{}_{t=0}={\vec{u}^{in}}.
\end{array}\right.
\end{equation}

Note that given
$\vec{u}=\overline{\vec{u}^{}} - \alpha^2 \Delta \overline{\vec{u}^{}}$ 
 the Poincar\'e inequality $ \displaystyle \|{\vec{u}} \|^{}_{L^{2}} \le {L}/{2\pi} \|\nabla \vec{u} \|_{L^{2}}$ immediately leads to
\begin{equation}
\label{isomorphisme}
\alpha^2 \|\overline{\vec{u}^{}}\|_{H^2} \le \|{\vec{u}^{}}\|_{L^2} \le  (\frac{L^2}{4\pi^2}+ \alpha^2)\|\overline{\vec{u}^{}}\|_{H^2}.
\end{equation}
In order to proof the ladder inequalities (\ref{inegalite d'echelle}) we need  first  to  show a regularity result for (\ref{Malpha ns}) or (\ref{Malphaexpanded ns}).
\begin{Proposition}
\label{regul}
Assume $\vec{f} \in {H}^{m-1}(\mathbb{T}_3)^3 $ and ${\vec{u}}^{in} \in {H}^{m}(\mathbb{T}_3)^3$, $ m \ge 1$, then the solution  $(\vec{u}^{},p)$ of (\ref{Malpha ns}) is such that 
\begin{align}
\vec{u} & \in L^{\infty}([0,T],{H}^{m}(\mathbb{T}_3)^3)\cap L^{2}([0,T],{H}^{m+1}(\mathbb{T}_3)^3),\\ 
p & \in  L^{2}([0,T],{H}^{m}(\mathbb{T}_3)^3) .
 \end{align}
\end{Proposition}
 The following Theorem  is a direct consequence of proposition \ref{regul}.

 \begin{Theorem}
\label{regularity}
 Assume  $\vec{f} \in C^{\infty}(\mathbb{T}_3)^{3} $ and ${\vec{u}}^{in} \in C^{\infty}(\mathbb{T}_3)^{3}$. Let $( \vec{u}^{},p^{})$ be the solution to problem (\ref{Malpha ns}). Then  the solution is $C^{\infty}$ in space and time.
\end{Theorem}
 The aim of this section is the proof of   proposition \ref{regul}.  We begin by recalling some known result for (\ref{Malpha ns}) or (\ref{Malphaexpanded ns}).
 \subsection{Known results}
Results in \cite{ILT05}  can be summarised as follows:
\begin{Theorem}
\label{existen}
 Assume $\vec{f} \in {L}^{2}(\mathbb{T}_3)^{3}$ and ${\vec{u}}^{in} \in
{H}^{1}(\mathbb{T}_3)^{3}$.  Then for any $T>0$, (\ref{Malpha ns}) has a unique distributional   solution $(\vec{u}^{ },{p}^{ }):=(\vec{u}^{\alpha },{p}^{\alpha })$  such that 
\begin{align}
{\vec{u}^{ }}  \in L_{}^{\infty}([0,T],{H}^{-1}(\mathbb{T}_3)^{3})\cap L^{2}([0,T],{L}^{2}(\mathbb{T}_3)^{3}),\\
\overline{\vec{u}^{}}  \in L_{}^{\infty}([0,T],{H}^{1}(\mathbb{T}_3)^{3})\cap L^{2}([0,T],{H}^{2}(\mathbb{T}_3)^{3}),
\end{align}
\begin{equation}
\label{thisone}
\begin{array}{lll}
 \displaystyle \|\overline{\vec{u}^{}}(t)\|_{L^2}^2+ \alpha^2 \| \overline{\vec{u}^{}}(t)\|_{H^1}^2 \le  \displaystyle (\|{\vec{u}^{in}}\|_{L^2}^2 
+ \alpha^2 \| {\vec{u}^{in}}\|_{H^1}^2)\exp{\left({-4\pi\nu t}/{L^2}\right)} \\ 
 \hskip 5cm \displaystyle + \frac{L^2}{4\pi^2\nu^2}\|\vec{f}\|_{H^{-1}}^2 \left(1 - \exp{\left({-4\pi\nu t}/{L^2}\right)}\right).
 \end{array}
\end{equation}
 Furthermore, if
${\vec{u}}^{in} \in
{H}^{2}(\mathbb{T}_3)^{3}$ then 
\begin{align}
& {\vec{u}^{ }}  \in L_{}^{\infty}([0,T],{L}^{2}(\mathbb{T}_3)^{3}),\\
& \overline{\vec{u}^{}}  \in L_{}^{\infty}([0,T],{H}^{2}(\mathbb{T}_3)^{3}),\\
\displaystyle & \|\overline{\vec{u}^{}}(t)\|_{H^1}^2 + \alpha^2 \| \overline{\vec{u}^{}}(t)\|_{H^2}^2 \le k(t) \label{titiimproved}.
\end{align}
Where  $k(t)$ verifies in particular:\\
(i) $k(t)$ is finite for all $t >0$.\\
(ii) $\displaystyle \limsup_{t\rightarrow \infty} k(t) < \infty.$
\end{Theorem}
\begin{Remark}
(1) The proof is based on the following energy inequality that is obtained by taking the inner product of (\ref{Malpha ns}) with $\overline{\vec{u}},$
\begin{equation}
\begin{array}{llcc}
 \label{energy equality  l2 bar ns alpha}
 \displaystyle \frac{1}{2} (
\frac{d}{dt}
\|\overline{\vec{u}}^{}\|_{{L}^{2}}^{2}+\alpha^{2} \frac{d}{dt}
\|\nabla \overline{\vec{u}}^{}\|_{{L}^{2}}^{2}) +\nu (\|\nabla
\overline{\vec{u}}^{}\|_{{L}^{2}}^{2}+\alpha^{2}\|\Delta
\overline{\vec{u}}^{}\|_{{L}^{2}}^{2} )\le \| \vec{f}\|_{{L}^{2}}\| \overline{\vec{u}}^{}\|_{{L}^{2}}.
\end{array}
\end{equation}
(2) Note that the pressure may be reconstructed from ${\vit}$ and $\overline{\vit}$ by solving the elliptic equation
$$\Delta p = \nabla \cdot ( ( {\vit} \cdot \nabla)\overline{\vit}).$$
One
concludes from the classical elliptic theory that $p \in  L^{1}([0,T],{L}^{2}(\mathbb{T}_3)^{3}).$
\end{Remark}
We recall that we can extract subsequences of solution that converge as $\alpha \rightarrow 0$ to a weak solution of the Navier Stokes equations.
The reader can look in  \cite{ILT05},  \cite{FDT02} and \cite{A09} for more details.  
\begin{Corollary}
\label{sobolev}
(1) We have ${\vit}  \in
\displaystyle  L^{2}([0,T],L^{2}(\mathbb{T}_3))$ and by Sobolev embending $\overline{\vit} \in
\displaystyle  L^{2}([0,T],L^{\infty}(\mathbb{T}_3))$. Thus there exists a constant $M(T):=M(\vec{u}^{in}, \vec{f}, \alpha, T) > 0$ such that $$ \displaystyle \int_{0}^{t}\|\overline{\vec{u}} \|^{2}_{L^{\infty}} \le \frac{1}{\alpha^2}\int_{0}^{t}\|{\vec{u}} \|^{2}_{L^{2}} \le  M(T)  \ \ \hbox{ for all }   t \in [0,T].$$
(2) We also observe by using (\ref{isomorphisme})  that there exists a constant $C(\alpha):=C(\alpha,L) >0$  such that 
\begin{align}
\displaystyle  \|{\vec{u}^{}}(t)\|_{L^2}^2 \le C(\alpha)  k(t) \hbox{ for all }   t >0.
\end{align}
 \end{Corollary}
\subsection{Regularity: Proof of proposition  \ref{regul} }
The proof of proposition  \ref{regul} is classical (see for example in \cite{LL06b}). In order to make the paper self-contained we will give a complete proof for this regularity result. The
proof is given in many steps.\\
\textbf{Step 1:}  we show that  ${\vec{u}} \in L^{\infty}([0,T],L^{2}(\mathbb{T}_3)^{3})\cap L^{2}([0,T],{H}^{1}(\mathbb{T}_3)^{3})$.\\
\textbf{Step 2:} we take $\partial_{t}\vec{u}$ as a test function in (\ref{Malpha ns}).\\
\textbf{Step 3:} We take  the $m-1$  derivative of    (\ref{Malpha ns}) then we take  $\partial_{t}\nabla^{m-1}\vec{u}$ as a test function and the result follows by induction.\\
  
\textbf{Step 1:}\\
We have the following  Lemma.
\begin{lemme}
\label{propp}
For  ${\vec{u}}^{in} \in
{L}^{2}(\mathbb{T}_3)^{3}$ and  $\vec{f} \in {H}^{-1}(\mathbb{T}_3)^{3}$, eqs. (\ref{Malpha ns}) have a unique solution $ (\vit,p)$ such that
\begin{align}
  &{\vec{u}} \in L^{\infty}([0,T],L^{2}(\mathbb{T}_3)^{3})\cap L^{2}([0,T],{H}^{1}(\mathbb{T}_3)^{3}).
  \end{align}
\end{lemme}
\textbf{Proof of Lemma \ref{propp}.}
We show formal a priori estimates for the solution established in  Theorem \ref{existen}. These estimates can be obtained rigorously using the Galerkin procedure.\\
We take the inner product of (\ref{Malpha ns}) with $\vec{u}$ to obtain
\begin{equation}
\begin{array}{llcc}
 \label{energy equality  l2 ns alpha}
 \displaystyle \frac{1}{2}
\frac{d}{dt}
\|{\vec{u}}(t,\textbf{x})\|_{{L}^{2}}^{2} +\nu \|\nabla
{\vec{u}}(t,\textbf{x})\|_{{L}^{2}}^{2}  \le \| \nabla^{-1}\vec{f}\|_{{L}^{2}}\| \nabla\vec{u}\|_{{L}^{2}}+ |(  (\vec{u}
\cdot \nabla)\overline{\vec{u}} , \vec{u} )|.
\end{array}
\end{equation}
Integration by parts and  Cauchy-Schwarz inequality yield to
 \begin{equation} |(  (\vec{u}
\cdot \nabla)\overline{\vec{u}} , \vec{u} )|\le  \|\vec{u}\otimes
\overline{\vec{u}} \|_{{L}^{2}} \| \nabla \vec{u}\|_{{L}^{2}}\end{equation}
and by Young's inequality, we obtain
 \begin{equation}
 \begin{array}{cc}
 \displaystyle \| \nabla^{-1}\vec{f}\|_{{L}^{2}}\| \nabla\vec{u}\|_{{L}^{2}}\le \frac{1}{\nu}\| \nabla^{-1}\vec{f}\|_{{L}^{2}}^{2}+\frac{\nu}{4}\| \nabla\vec{u}\|_{{L}^{2}}^{2},\\
 \displaystyle |(  (\vec{u}
\cdot \nabla)\overline{\vec{u}} , \vec{u})|\le \frac{1}{\nu} \|\vec{u}\otimes
\overline{\vec{u}} \|_{{L}^{2}}^{2}+ \frac{\nu}{4}\| \nabla \vec{u}\|_{{L}^{2}}^{2}.
\end{array}
\end{equation}
From the above inequalities we get
 \begin{equation}
 \label{jfkjfkjkdv}\begin{array}{llll}
  \displaystyle \frac{d}{dt}
\|{\vec{u}}(t,\textbf{x})\|_{{L}^{2}}^{2} +\nu \|\nabla
{\vec{u}}(t,\textbf{x})\|_{{L}^{2}}^{2}  &\le& \displaystyle
 \frac{2}{\nu}\| \nabla^{-1}\vec{f}\|_{{L}^{2}}^{2}
+ \frac{2}{\nu} \|\vec{u}
\overline{\vec{u}} \|_{{L}^{2}}^{2}\\
&\le& \displaystyle \frac{2}{\nu}\| \nabla^{-1}\vec{f}\|_{{L}^{2}}^{2}
+ \frac{2}{\nu}\frac{1}{\alpha^2} \|\vec{u}\|_{{L}^{2}}^{4},
\end{array}
\end{equation}
where we have used in the last step that 
\begin{equation}
\|\overline{\vec{u}} \|_{{L}^{\infty}}^{2} \le \frac{1}{\alpha^2} \|\vec{u}\|_{{L}^{2}}^{2}.
\end{equation}
This implies that
  \begin{equation}
 \label{jfkjfkjkdv12154}\begin{array}{llll}
  \displaystyle \frac{d}{dt}(1+
\|{\vec{u}}(t,\textbf{x})\|_{{L}^{2}}^{2} ) \le 
C_1 (1+ \|\vec{u}(t,\textbf{x})\|_{{L}^{2}}^{2})^2,
\end{array}
\end{equation}
where $C_1 = \max (\frac{2}{\nu}\frac{1}{\alpha^2},\displaystyle \frac{2}{\nu}\| \nabla^{-1}\vec{f}\|_{{L}^{2}}^{2} )$,
and by  Gronwall's Lemma, since $\|{\vec{u}} \|^{2}_{L^{2}}\in L^{1}([0,T]) $  (Corollary \ref{sobolev}) we conclude that
$$ 1+\|{\vec{u}}(t,\textbf{x})\|_{L^{\infty}([0,T],{L}^{2})}^{2} \le K_{1}(T),$$ where $ K_{1}(T):=K_{1}(T, \vec{u}^{in}, \vec{f})$ is given by  $$  K_{1}(T)=(1+\|{\vec{u}}^{in}\|_{{L}^{2}}^{2} )\exp{\left({C_1\int_{0}^{T} (1 + \|{\vec{u}}(s) \|^{2}_{L^{2}}) ds}\right)}. 
$$
Furthermore, for every $T>0$ we have from (\ref{jfkjfkjkdv}),
\begin{equation}
\|{\vec{u}}(T,\textbf{x})\|_{{L}^{2}}^{2} +\nu \int_{0}^{T}\|\nabla
{\vec{u}}(t,\textbf{x})\|_{{L}^{2}}^{2} dt \le \|{\vec{u}}^{in}\|_{{L}^{2}}^{2}
 +\frac{2}{\nu}\| \nabla^{-1}\vec{f}\|_{{L}^{2}}^{2}T
+ \frac{2}{\nu} K_{1}M.
\end{equation}
Thus $\vec{u} \in  L^{2}([0,T],{H}^{1}(\mathbb{T}_3)^{3})$ for all $T>0$.\\

\textbf{Step 2:}\\
With the same assumption in the inital data as in Theorem \ref{existen} , we can find the following result:
\begin{lemme}
\label{existence}
 Assume $\vec{f} \in {L}^{2}(\mathbb{T}_3)^{3}$ and ${\vec{u}}^{in} \in
{H}^{1}(\mathbb{T}_3)^{3}$ . Then for any $T>0$,  eqs. (\ref{Malpha ns})   have a unique regular  solution $(\vit,p)$ such that  
\begin{align}
 &{\vec{u}} \in C([0,T],{H}^{1}(\mathbb{T}_3)^{3})\cap L^{2}([0,T],{H}^{2}(\mathbb{T}_3)^{3}),\\
 &\displaystyle \frac{d {\vec{u}}}{dt } \in L^{2}([0,T],{L}^{2}(\mathbb{T}_3)^{3}),\\
 &p \in  L^{2}([0,T],{H}^{1}(\mathbb{T}_3)^{3}).
  \end{align}
\end{lemme}
\textbf{Proof of Lemma \ref{existence}}
\ It is easily checked that since $\vec{u}
\in L^{\infty}([0,T],{L}^{2}(\mathbb{T}_3)^{3}) \cap
 L^{2}([0,T],{H}^{1}(\mathbb{T}_3)^{3})$, then $\overline{\vit} \in L^{\infty}([0,T],{H}^{2}(\mathbb{T}_3)^{3}) \cap L^{2}([0,T],{H}^{3}(\mathbb{T}_3)^{3})$. Consequently, by
Sobolev injection Theorem, we deduce that $\overline{\vit} \in
\displaystyle L^{\infty}([0,T],L^{\infty}(\mathbb{T}_3)^{3})$ and $\nabla \overline{\vit} \in  L^{2}([0,T],L^{\infty}(\mathbb{T}_3)^{3})$.  \\
Therefore,
\begin{equation}
\label{l2}
 ({\vit}\cdot \nabla)\overline{\vit}^{} \in
L^{2}([0,T],L^{2}(\mathbb{T}_3)^{3}).\end{equation}
Now,
for fixed t, we can   take $\partial_{t}\vec{u}$ as a test
function in (\ref{Malpha ns}) and the procedure is the same as the one in \cite{RL}. Note that the proof given in \cite{RL} is formal and can be obtained rigorously by using Galerkin method combined with (\ref{l2}).

Once we obtain
that $ \displaystyle \vec{u} \in L^{\infty}([0,T],{H}^{1}(\mathbb{T}_3)^{3}) \cap L^{2}([0,T],{H}^{2}(\mathbb{T}_3)^{3}) \cap H^{1}([0,T],{L}^{2}(\mathbb{T}_3)^{3})$ and $p \in L^{2}([0,T],{H}^{1}(\mathbb{T}_3)^{3})$. Interpolating between $ L^{2}([0,T],{H}^{2}(\mathbb{T}_3)^{3})$ and $H^{1}([0,T],{L}^{2}(\mathbb{T}_3)^{3})$ yields to $\vec{u} \in C([0,T],{H}^{1}(\mathbb{T}_3)^{3})$.\\

\textbf{Step 3:}\\
\ We proceed by induction. The case $m = 1$ follows from Lemma \ref{existence}.\\
Assume that for any $k=1,...,m-1$, if   $\vec{f} \in {H}^{k-1}(\mathbb{T}_3)^{3} $ and ${\vec{u}}^{in} \in {H}^{k}(\mathbb{T}_3)^{3}$ then  $\vec{u}  \in L^{\infty}([0,T],{H}^{k}(\mathbb{T}_3)^{3})\cap L^{2}([0,T],{H}^{k+1}(\mathbb{T}_3)^{3})$ holds.\\
It remains to prove  that when  $k=m$, $\vec{f} \in {H}^{m-1}(\mathbb{T}_3)^{3} $ and ${\vec{u}}^{in} \in {H}^{m}(\mathbb{T}_3)^{3}$ that  $\vec{u}  \in L^{\infty}([0,T],{H}^{m}(\mathbb{T}_3)^{3})\cap L^{2}([0,T],{H}^{m+1}(\mathbb{T}_3)^{3})$.\\
 It is easily checked that for $\displaystyle\vec{u}
\in L^{\infty}([0,T],{H}^{k}(\mathbb{T}_3)^{3}) \cap
 L^{2}([0,T],{H}^{k+1}(\mathbb{T}_3)^{3})$, \\
 $\displaystyle\overline{\vit} \in L^{\infty}([0,T],{H}^{k+2}(\mathbb{T}_3)^{3}) \cap L^{2}([0,T],{H}^{k+3}(\mathbb{T}_3)^{3})$. Consequently, by
Sobolev injection Theorem, we deduce that $\nabla ^{k}\overline{\vit} \in
\displaystyle L^{\infty}([0,T],L^{\infty}(\mathbb{T}_3)^{3}),$ and $\nabla ^{k+1}\overline{\vit} \in  L^{2}([0,T],L^{\infty}(\mathbb{T}_3)^{3})$.  \\
By taking the $m-1$ derivative of (\ref{Malpha ns}) we get in the sense of the distributions that 
\begin{equation}
\label{nabla Malpha ns}
 \left\{
\begin{array} {llll} \displaystyle
 \frac{\partial  \nabla^{m-1} {\vec{u}^{}} }{\partial t}+  \nabla^{m-1} \left(({\vit}\cdot \nabla)\overline{\vit}^{}\right)
 - \nu  \nabla^{m-1} \Delta {\vec{u}^{}}  + \nabla^{m-1}\nabla p^{} =
{ \nabla^{m-1} \vec{f}^{}},\\
 \nabla \cdot  \nabla^{m-1}{\vit}=0,\\
 \nabla^{m-1}{\vec{u}^{}}_{t=0}={ \nabla^{m-1}\vec{u}^{in}}.
\end{array}\right.
\end{equation}
where boundary conditions remain periodic and still with zero mean and the initial
condition with zero divergence and mean.
 
Therefore, after using Leibniz Formula
\begin{equation}
\label{l212123}
 \nabla^{m-1} \left(({\vit}\cdot \nabla)\overline{\vit}^{}\right)= \sum_{k=0}^{m-1} C_{m-1}^{k}\nabla^{k}\vec{u}\nabla^{m-k} \overline{\vec{u}}^{},
\end{equation}
since $$\nabla^{k}\vec{u}
\in L^{\infty}([0,T],{L}^{2}(\mathbb{T}_3)^{3})$$ and $$\nabla ^{k+1}\overline{\vit} \in  L^{2}([0,T],L^{\infty}(\mathbb{T}_3)^{3}),$$ for any $k=1,...,m-1.$\\
We deduce that
\begin{equation}
\label{l2primeprime}
 \nabla^{m-1} \left(({\vit}\cdot \nabla)\overline{\vit}^{}\right) \in
L^{2}([0,T],L^{2}(\mathbb{T}_3)^{3}).
\end{equation}

Now,
for fixed t, we can   take  $\partial_{t}\nabla^{m-1}\vec{u}$ as a test
function in (\ref{nabla Malpha ns}) and the procedure is the same as the one in \cite{RL}.   One obtains
that $ \displaystyle \vec{u} \in L^{\infty}([0,T],{H}^{m}(\mathbb{T}_3)^{3}) \cap L^{2}([0,T],{H}^{m+1}(\mathbb{T}_3)^{3})$ and $p \in L^{2}([0,T],{H}^{m}(\mathbb{T}_3)^{3})$.
This finishes the proof of proposition \ref{regul}.

\section{ Ladder Inequalities:  Proof of theorem \ref{thhhh}.}
The first step in the proof of theorem \ref{thhhh}, which has been expressed in section 1, is the  energy inequality (\ref{energy equality  l2 bar ns alpha}) that  corresponding to the case $N=0$ of  (\ref{inegalite d'echelle}).
Having showing  in the above section the regularity result for (\ref{Malpha ns}). We can take the $N$ derivative of (\ref{Malpha ns}),
 we get in the sense of the distributions  that for all $N \ge 1$,
\begin{equation}
\label{nablaprime Malpha ns}
 \left\{
\begin{array} {llll} \displaystyle
 \frac{\partial  \nabla^{N} {\vec{u}^{}} }{\partial t}+  \nabla^{N} \left(({\vit}\cdot \nabla)\overline{\vit}^{}\right)
 - \nu  \nabla^{N} \Delta {\vec{u}^{}}  + \nabla^{N}\nabla p^{} =
{ \nabla^{N} \vec{f}^{}},\\
 \nabla \cdot  \nabla^{N}{\vit}=0,\\
 \nabla^{N}{\vec{u}^{}}_{t=0}={ \nabla^{N}\vec{u}^{in}}.
\end{array}\right.
\end{equation}
where boundary conditions remain periodic and still with zero mean and the initial
condition with zero divergence and mean.
%
%
Taking $ \nabla^N \bu $ as test function in  (\ref{nablaprime Malpha ns}),
we can write that
$$ \frac{1}{2}\frac{d}{dt}\int_{\mathbb{T}_3}\left|\nabla^{N}\overline{\vec{u}}\right|^{2} d\textbf{x}+ \alpha^{2} \frac{1}{2}\frac{d}{dt}\int_{\mathbb{T}_3}\left|\nabla^{N+1}\overline{\vec{u}}\right|^{2} d\textbf{x}=\nu \int_{\mathbb{T}_3}\nabla^{N}\overline{\vec{u}}   \nabla^{N} \Delta \overline{\vec{u}}d\textbf{x}-\nu\alpha^{2}\int_{\mathbb{T}_3}\nabla^{N}\overline{\vec{u}} \nabla^{N}\Delta   \Delta\overline{\vec{u}}d\textbf{x} $$
$$+\int_{\mathbb{T}_3}\nabla^{N}\overline{\vec{u}} \nabla^{N}( (\overline{\vec{u}} \cdot \nabla ) \overline{\vec{u}})d\textbf{x}-\alpha^{2}
\int_{\mathbb{T}_3}\nabla^{N}\overline{\vec{u}}
\nabla^{N}((\Delta\overline{\vec{u}} \cdot
\nabla)\overline{\vec{u}})d\textbf{x}
+\int_{\mathbb{T}_3}\nabla^{N}\overline{\vec{u}}
\nabla^{N}{\vec{f}}d\textbf{x}.$$
Where the pressure term vanishes as $\nabla \cdot \nabla^N \bu=0$. 


Using the defintion of $\overline{H_N}$ in (\ref{propriete10}) we obtain 
 \begin{equation}
\label{....................5...........}
 \begin{array}{cc}
 \displaystyle \frac{1}{2}(\frac{d}{dt}{\overline{H_{N}}}+\displaystyle \alpha^{2}\frac{d}{dt}{\overline{H_{N+1}}}) \le \displaystyle -\nu (\overline{H_{N+1}}+\displaystyle \alpha^{2}\overline{H_{N+2}})
+|\displaystyle \int_{\mathbb{T}_3}\nabla^{N}\overline{\vec{u}}
\nabla^{N}((\overline{\vec{u}} \cdot \nabla)
\overline{\vec{u}})d\textbf{x}|\\
\hskip 5 cm +\displaystyle \alpha^{2}|\int_{\mathbb{T}_3}\nabla^{N+1}\overline{\vec{u}}
\nabla^{N-1}((\Delta\overline{\vec{u}} \cdot
\nabla)\overline{\vec{u}})d\textbf{x}|+\displaystyle |\int_{\mathbb{T}_3}\nabla^{N}\overline{\vec{u}}
\nabla^{N}{\vec{f}}d\textbf{x}|.
\end{array}
\end{equation}
Where we have integrated by parts in the Laplacien terms.\\
The central terms are 
 \begin{equation}
 \hbox{NL}_1=\displaystyle |\int_{\mathbb{T}_3}\nabla^{N}\overline{\vec{u}} \nabla^{N}((\overline{\vec{u}} \cdot \nabla ) \overline{\vec{u}})d\textbf{x}|
 \end{equation}
 and 
 \begin{equation}
 \hbox{NL}_2=\displaystyle \displaystyle \alpha^{2}|\int_{\mathbb{T}_3}\nabla^{N+1}\overline{\vec{u}}
\nabla^{N-1}((\Delta\overline{\vec{u}} \cdot
\nabla)\overline{\vec{u}})d\textbf{x}|
 \end{equation}
These two terms   $ \hbox{NL}_1 $ and $ \hbox{NL}_2$  can be bounded using the following  Gagliardo-Nirenberg interpolation inequality \cite{DG95}: 
 \begin{lemme}
 The Gagliardo-Nirenberg interpolation inequality is:\\ 
  For
 $\displaystyle 1 \le q,r \le \infty $, $j$ and  $m$ such that
 $\displaystyle 0 \le j < m $ we have
\begin{equation}
\label{gag}
 \left\|\nabla^{j} \vec{v}\right\|_{p} \le  C \left\|\nabla^{m} \vec{v} \right\|_{r}^{a}   \left\| \vec{v} \right\|_{q}^{1-a}
\end{equation}
where $$ \frac{1}{p}=\frac{j}{d}+a \left( \frac{1}{r}-\frac{m}{d}
\right) +\frac{1-a}{q}$$
for   $\displaystyle \frac{j}{m} \le a < 1$ and $\displaystyle a = \frac{j}{m}$  if  $\displaystyle m -j - \frac{d}{r} \in N^{\ast} $ .\\
 \end{lemme}
The first nonlinear term  $\hbox{NL}_1 $ 
 is estimated with the Gagliardo-Nirenberg inequality \cite{DG95} by $\displaystyle c_{N}   \left\|\nabla\overline{\vec{u}}\right\|_{\infty} \overline{H_{N}}$, where $c_{0}=0$ and $c_{N}  \le c 2^{N}$. Indeed, the  nonlinear first term 
  $\hbox{NL}_1 $ 
 is found to satisfy   $$\hbox{NL}_1 =\left|\int_{\mathbb{T}_3}\nabla^{N}\overline{\vec{u}} \nabla^{N}((\overline{\vec{u}} \cdot \nabla ) \overline{\vec{u}})d\textbf{x} \right|
\le  2^{N} \overline{H_{N}}^{1/2}\sum_{l=1}^{N} \left\|\nabla^{l}\overline{\vec{u}}^{}\right\|_{L^{p}} \left\|\nabla^{N+1-l}\overline{\vec{u}}^{} \right\|_{L^{q}},
$$
where $p$ and $q$ satisfy $1/p+1/q=1/2$ according to the {H\"{o}lder} inequality.
  We use now the two Gagliardo-Nirenberg
inequalities
$$
 \left\|\nabla^{l}\overline{\vit}\right\|_{L^{p}} \le c_{1} \left\|\nabla^{N}\overline{\vit}\right\|_{L^{2}}^{a}\left\|\nabla \overline{\vit}\right\|_{\infty}^{1-a},
$$
$$
\left\|\nabla^{N+1-l}\overline{\vec{u}}^{} \right\|_{L^{q}} \le
c_{2} \left\|\nabla^{N}\overline{\vec{u}}^{}\right\|_{L^{2}}^{b}\left\|\nabla\overline{\vec{u}}^{}\right\|_{\infty}^{1-b}.
$$
Where $a$ anb $b$ must satisfy
$$
 \frac{1}{p}=\frac{l-1}{3}+a \left( \frac{1}{2}-\frac{N-1}{3} \right),
$$
$$\frac{1}{q}=\frac{N-l}{3}+b \left( \frac{1}{2}-\frac{N-1}{3} \right).$$
Since $1/p+1/q=1/2$, we deduce  $a + b =1.$
Thus we obtain
 \begin{equation}
 \left|\int_{\mathbb{T}_3}\nabla^{N}\overline{\vec{u}} \nabla^{N}((\overline{\vec{u}} \cdot \nabla  )\overline{\vec{u}})d\textbf{x} \right|
\le  c_{N}   \left\|\nabla\overline{\vec{u}}\right\|_{\infty} \overline{H_{N}} .
\end{equation}

In the same way, we can estimate the nonlinear second term with
Gagliardo-Nirenberg inequality  in order to have
 \begin{equation}
\label{.........d...........5...........}
 \alpha^{2}|\int_{\mathbb{T}_3}\nabla^{N+1}\overline{\vec{u}} \nabla^{N-1}((\Delta\overline{\vec{u}} \cdot \nabla)\overline{\vec{u}})d\textbf{x}|\le c^{'}_{N} \alpha^{2}  \left\|\nabla\overline{\vec{u}}\right\|_{\infty} \overline{H_{N+1}},
\end{equation}
 where  $c^{'}_{N} \le  c2^{N}.$\\
The result  (\ref{inegalite d'echelle}) then follows.

\section{Estimates in terms of Reynolds number: Proof of Theorem \ref{lengththoerem}}
\subsection{Proof of inequality (\ref{lengthimproved})  }

We begin by forming the combination 
$$ \overline{F_N} = \overline{H_N} + \tau {\Phi_{N}} $$
where the quantity $ \tau$ is defined by
$$ \tau = \ell^2 \nu^{-1} (Gr \ln{Gr} )^{-1/2}.$$
We define the combination 
$$J_N = \overline{F_N} + 2 \alpha^2 \overline{F_{N+1}}$$
The following result is a consequence of Theorem \ref{thhhh} and it is proof follows closely to that for the Navier Stokes-alpha model  in \cite{GH06} and we will not repeat it.  
\begin{Theorem}
\label{th2} As $Gr\rightarrow \infty$, for $N \ge 1$, $1 \le  p \le N$ the unique  solution to eqs. (\ref{Malpha ns}) satisfies
\begin{equation}
\label{3.5 3.6}
\begin{array}{ll}
 \displaystyle \frac{1}{2}\frac{d}{dt}{{J_{N}}} \le \displaystyle  -\nu \frac{J_N^{1+\frac{1}{p}}}{J_{N-p}^{\frac{1}{p}}}
\displaystyle + C_{N,\alpha}   \left\|\nabla\overline{\vec{u}}\right\|_{\infty} {J_{N}}  + C \nu \ell^{-2}Re (\ln{Re})J_N
 \end{array}
\end{equation}
and for $N=0,$
\begin{equation}
\label{3.5 3.6 4.5}
\begin{array}{ll}
 \displaystyle \frac{1}{2}\frac{d}{dt}{{J_{0}}} \le \displaystyle  -\nu J_1 +
\displaystyle   C \nu \ell^{-2}Re (\ln{Re})J_0
 \end{array}
\end{equation}
\end{Theorem}
When $ \alpha \rightarrow 0 $, $J_N $ tends to $F_N= {H_N} + \tau {\Phi_{N}}$, and the result of theorem \ref{th2} is consistent with the result achived to Navier Stokes equations in \cite{DG95}.\\
To obtain  length scales estimates  let us 
 define the quantities 
 $$ \kappa_{N,r}=\left(  \frac{J_N}{J_r} \right)^{\frac{1}{2(N-r)}}  $$ 
 
 In the $\alpha \rightarrow 0 $ limit, the $\kappa_{N,0}$ behaves as the $2N$th moment of the energy spectrum. \\ 
 The aim of this subsection is to find  a estimate for the length sacles associated with higher order moments solution $ \kappa_{N,0}$ $(N \ge  2 )$. To this end,   we find first  upper bounds for $\left\langle \kappa_{N,r}^2 \right\rangle $, $\left\langle \kappa_{1,0}^2 \right\rangle $ and $\left\langle \|\nabla \overline{\vec{u}}\|_{\infty} \right\rangle .$ Then we use the following identity 
 
  \begin{equation}
\begin{array}{llllll}
  \kappa_{N,0}^2 
&=&     \kappa_{N,1}^{2(N-1)/N} \kappa_{1,0}^{2/N}  
\end{array}
\end{equation}
in order to deduce the result. 

 (a) The First two bounds are  obtained by dividing by $J_N$ in Theorem \ref{th2} and time averaging to obtain 
 \begin{equation}
 \label{meme}
 \left\langle \kappa_{N,r}^2 \right\rangle \le   C_{N,\alpha} \nu^{-1} \left\langle  \left\|\nabla\overline{\vec{u}}\right\|_{\infty}  \right\rangle   + C\ell^{-2}Re (\ln{Re})
 \end{equation}
 and 
  \begin{equation}
  \label{k 1 0}
 \left\langle \kappa_{1,0}^2 \right\rangle \le  C\ell^{-2}Re (\ln{Re}).
 \end{equation}
 
 \begin{Remark}
 Note that the bound on $\left\langle \| \overline{\vec{u}}\|_{\infty}^2 \right\rangle $ is found to satisfies (see in \cite{GH06} for more details),
   \begin{equation}
   \label{u infini}
\left\langle \| \overline{\vec{u}}\|_{\infty}^2 \right\rangle \le C \ell^{-2}\nu^2 V_{\alpha}Re^{11/4}
 \end{equation}
 Where $$  V_{\alpha}:= \left(\frac{L}{(\ell \alpha)^{1/2}}^3\right).$$ 
 \end{Remark}
 
 (b) It is also possible to estimate  $\left\langle \|\nabla \overline{\vec{u}}\|_{\infty} \right\rangle $: return to the eqs. (\ref{Malpha ns}) and take a different way. We take $\vec{u}=-\alpha^2 \Delta \overline{\vec{u}} + \overline{\vec{u}} $ as test function,  then integration by parts (see Lemma \ref{propp}), using (\ref{titiimproved}) and time averaging,  we obtain 
\begin{equation}
\begin{array}{llllll}
\nu \left\langle \overline{H_1}+ 2 \alpha^2\overline{H_2} + \alpha^4 \overline{H_3}   \right\rangle 
&\le&  C \left\langle  \| \nabla \overline{\vec{u}}\|_{L^2} \| \Delta {\overline{\vec{u}}}\|_{L^2}^2 \right\rangle + (1+ \alpha^2 \ell^{-2}) \left\langle \overline{{H_{0}}}^{1/2} \textbf{$\Phi_{0}$}^{1/2}\right\rangle
\\
&\le&  C \left\langle  \| \Delta \overline{\vec{u}}\|_{L^2}^2\right\rangle^{} \| \overline{\vec{u}}\|_{L^{\infty}([0,T], H^1)}  + (1+ \alpha^2 \ell^{-2}) \left\langle \overline{{H_{0}}}^{1/2} \textbf{$\Phi_{0}$}^{1/2}\right\rangle\\
&\le&  C  \alpha^{-2}  \nu^2  \ell^{-4} L^3  Re^3  Gr^{2}                                   + C (1+ \alpha^2 \ell^{-2}) \nu^3  \ell^{-4} L^3  Re^3.
\end{array}
\end{equation}
 Thus we can right
 \begin{equation}
 \label{H 3 bar}
\begin{array}{llllll}
 \left\langle  \overline{H_3}   \right\rangle 
&\le&  C (\alpha, \nu , \ell, L) Re^7.
\end{array}
\end{equation}
This can be used to find the estimate for $\left\langle \|\nabla \overline{\vec{u}}\|_{\infty} \right\rangle $. In fact, Agmon's inequality \cite{FMRT01}
$$\displaystyle \|\vec{u}\|_{\infty} \le  \|\vec{u}\|_{H^1}^{1/2}  \|\vec{u}\|_{H^2}^{1/2} $$
  says that 
 \begin{equation}
 \label{nabla u infini}
\begin{array}{llllll}
 \left\langle   \|\nabla \overline{\vec{u}}\|_{\infty}  \right\rangle 
&\le& \left\langle  \overline{H_3}  \right\rangle^{1/4} \left\langle  \overline{H_2}  \right\rangle^{1/4}\\
&\le& C (\alpha, \nu , \ell, L) Re^{5/2}.
\end{array}
\end{equation}
 
(c) Thus we obtain from (\ref{meme}) and  (\ref{nabla u infini}) that 
   \begin{equation}
   \label{k n r}
\begin{array}{llllll}
\ell^{2} \left\langle   \kappa_{N,r}^2 \right\rangle 
&\le& C (\alpha, \nu , \ell, L) Re^{5/2}+CRe (\ln{Re}).
\end{array}
\end{equation}
In particular, for $r=0$ 
 \begin{equation}
   \label{k n 0}
\begin{array}{llllll}
\ell^{2} \left\langle   \kappa_{N,0}^2 \right\rangle 
&\le& C (\alpha, \nu , \ell, L) Re^{5/2}+CRe (\ln{Re}).
\end{array}
\end{equation}
By choosing $r=1$ we can then get an improvement for  $ \left\langle   \kappa_{N,0}^2 \right\rangle  $   by writting
   \begin{equation}
\begin{array}{llllll}
 \left\langle   \kappa_{N,0}^2 \right\rangle 
&=&  \left\langle   \kappa_{N,1}^{2(N-1)/N} \kappa_{1,0}^{2/N} \right\rangle \\
&\le&  \left\langle   \kappa_{N,1}^{2}  \right\rangle^{(N-1)/N} \left\langle \kappa_{1,0}^{2} \right\rangle^{1/N}, 
\end{array}
\end{equation}
and then using the above estimates for $\left\langle   \kappa_{N,1}^{2}  \right\rangle$ and $ \left\langle \kappa_{1,0}^{2} \right\rangle$ which give for $N \ge 2$, 
    \begin{equation}
    \label{N DEPEND}
\begin{array}{llllll}
\ell^{2} \left\langle   \kappa_{N,0}^2 \right\rangle 
&\le&  C(\alpha,\nu , \ell, L)^{(N-1)/N} Re^{5/2-3/2N} (\ln{Re})^{1/N}+ C Re \ln{Re}.
\end{array}
\end{equation}
Note that when $N=1$ we return to $  
\ell^{2} \left\langle   \kappa_{1,0}^2 \right\rangle 
\le   C Re \ln{Re}.
 $

\subsection{Proof of inequality (\ref{improvedbound kolomogrov})}

Let us come back to relation (\ref{inegalite d'echelle}) ,  when $N=0$,  
 we get the energy inequality (\ref{energy equality  l2 bar ns alpha}) 
\begin{equation}
\label{here}
\begin{array}{llllll}
\displaystyle \frac{d}{dt}({\overline{H_{0}}}+ \alpha^{2}{\overline{H_{1}}}  ) \le   
  -\displaystyle\nu ({\overline{H_{1}}}+ \alpha^{2} {\overline{H_{2}}} )
 + \overline{{H_{0}}}^{1/2} \textbf{$\Phi_{0}$}^{1/2}.
\end{array}
\end{equation}
Poincar\'e inequality together with Cauchy Schwarz, Young and Gronwall inequalities in (\ref{here}) imply that $ \displaystyle \overline{H_{0}}+ \alpha^{2}\overline{H_{1}} $ is uniformly bounded in time according to (\ref{thisone}).  
Time averaging, using the fact that the time average of the time derivative in (\ref{here}) vanishes,   we obtain 
\begin{equation}
\begin{array}{llllll}
\displaystyle   \left\langle  {\overline{H_{1}}}+ \alpha^{2}{\overline{H_{2}}} \right\rangle
&\le&   \left\langle \overline{{H_{0}}}^{1/2} \textbf{$\Phi_{0}$}^{1/2}\right\rangle\\
&\le& c\nu^2  \ell^{-4} L^3 Re^3.
\end{array}
\end{equation}
Thus
\begin{equation}
\label{H 1 prime bar}
\begin{array}{llllll}    \left\langle    {\overline{H_{1}}} \right\rangle 
&\le& c   \nu^2  \ell^{-4} L^3  Re^3.
\end{array}
\end{equation}
and 
\begin{equation}
\label{H 2 bar}
\begin{array}{llllll}    \left\langle    {\overline{H_{2}}} \right\rangle 
&\le& c \alpha^{-2}  \nu^2  \ell^{-4} L^3  Re^3.
\end{array}
\end{equation}

From (\ref{H 2 bar}) and by using the following interpolation inequality 
\begin{equation}
\begin{array}{llllll}   \displaystyle \overline{H_N} \le \displaystyle  \overline{H_{N-s}}^{\frac{r}{r+s}}  \displaystyle \overline{H_{N+r}}^{\frac{s}{r+s}},
\end{array}
\end{equation}
that is 
\begin{equation}
\begin{array}{llllll}   \displaystyle \overline{H_1} \le \displaystyle  \overline{H_{0}}^{\frac{1}{2}}  \displaystyle \overline{H_{2}}^{\frac{1}{2}},
\end{array}
\end{equation}
we can improve   (\ref{H 1 prime bar}) in order to obtain
\begin{equation}
\label{H 1 bar}
\begin{array}{llllll}   \displaystyle  \left\langle \overline{H_1}  \right\rangle \le \displaystyle  \left\langle \overline{H_{0}}\right\rangle^{\frac{1}{2}}  \displaystyle \left\langle \overline{H_{2}}\right\rangle^{\frac{1}{2}} \le  
\displaystyle c \alpha^{-1}  \nu^2  \ell^{-3} L^3  Re^{5/2}.
\end{array}
\end{equation}
This improve the Navier Stokes result (\ref{bound kolomogrov}) for the inverse Kolomogorov lentgh    
to 
\begin{equation}
\label{l lambda -1}
\ell \lambda_k^{-1} \le c (\frac{\ell}{\alpha})^{1/4} Re^{5/8} .
\end{equation}
We also deduce that  the energy dissipation rate $ \overline{\epsilon}=\nu \left\langle  \|\nabla \overline{\vec{u}}\|_{L^2}^2\right\rangle L^{-3} $ is also bounded by $Re^{5/2}$ but all the improved estimates blow up when $\alpha$ tends to zero.\\
\subsection{Proof of inequality  (\ref{attracteurimproved}) }
The authors in \cite{ILT05} showed the existence of a global attractor $ \mathcal{A}$  to this model and they made   estimates of the fractal dimension of this  attractor. 
The sharp estimate found  in \cite{ILT05} for the fractal dimension of $ \mathcal{A}$  expressed in terms of Grashof number $Gr$ is
\begin{equation}
\label{grashofattractor dimension}
d_{F, ML-\alpha}(\mathcal{A}) \le c \left( \frac{2 \pi}{L \alpha^2  \gamma}  \right)^{3/4} Gr^{3/2}.
\end{equation}
Where $$\displaystyle \frac{1}{\gamma}=\min{(1,\frac{2 \pi}{\alpha^2  L})}$$  
 This however can be improved by noting that their estimate  depends upon $\left\langle  {\overline{H_{1}}}+ \alpha^{2}{\overline{H_{2}}} \right\rangle$  whose upper bound is $Re^3$ not $Gr^2 \le c Re^4.$ 
With this improvement it is found that the estimate of $d_{F, ML-\alpha}(\mathcal{A})$ in \cite{ILT05} convert to 
\begin{equation}
\label{attractor dimension}
d_{F, ML-\alpha}(\mathcal{A}) \le c \left( \frac{L^{3/2} (2\pi)^{3/2}\ell^{-4}}{\alpha^2 }  \right)^{3/4} Re^{9/4}.
\end{equation}
In term of degrees of fredoom, this result says that $ Re^{3/4}\times Re^{3/4}\times Re^{3/4}$ resolution grid points are needed.
\\
 
\textbf{Acknowledgement}\\
\ \     \ \  I would like to thank Professor Roger Lewandowski and Professor J.D. Gibbon for helpful remarks and suggestions.

\end{document}